\documentclass{article}
\usepackage{amssymb}

\usepackage{amsmath}

\newtheorem{theorem}{Theorem}

\newtheorem{definition}[theorem]{Definition}

\newtheorem{lemma}[theorem]{Lemma}

\begin{document}

\title{Bismut's way of the Malliavin Calculus for non-markovian semigroups: an introduction.}
\author{R\'emi L\'eandre   \\
%EndAName
Laboratoire  de Math\'ematiques, \\ Universi\'e de Bourgogne-Franche-Comte,
\\ 25030, Besan\c con, FRANCE.
\\email: remi.leandre@univ-fcomte.fr
%\\
% \\Habib Ouerdiane\\ D\'epartement de Math\'ematiques. Facult\'e
%des Sciences de Tunis.\\ Universit\'e de Tunis El Manar. 1060.
%Tunis. TUNISIA
} \maketitle

\begin{abstract} We give a review of our recent works related to the Malliavin Calculus of Bismut type for non markovian generator. Part IV is new and relates the Malliavin Calculus and the general theory of elliptic pseudo-differential operators.\end{abstract}

{\section{Introduction}
Le $M$ be a compact Riemannian manifold endowed with its natural Riemannian measure $dx$ ($x$ is the generic element of $M$).  In local coordinates, we can think of the linear space $\mathbb{R}^d$ endowed with the metric $g_{i,j}(x) dx^i \otimes dx^j$  where $x \rightarrow (g_{.,.}(x))$ is a smooth function from $\mathbb{R}^d$ into the space of symmetric strictly positive matix. The Riemannian measure associated is
\begin{equation}dx = det(g_{.,.})^{-1/2}dx^1..\otimes dx^d\end{equation}
We consider a {\bf{ linear}} symmetric positive  operator densely defined on $L^2(dx)$  acting on a space which separates the point on $M$.  This means if $f$ 	and $g$ belong to this space,
\begin{equation}\int_Mg(x)Lh(x)dx = \int_Mh(x)Lg(x)dx\end{equation}
\begin{equation}\int_Mh(x)Lh(x)dx\geq 0\end{equation}
It has by abstract theory a selfadjoint extension on $L^2(dx)$, which generates a contraction semi-group $P_t$ on $L^2(dx)$ which solves the heat equation for $t>0$ 
\begin{equation}{\partial\over \partial t}P_th= -LP_th\end{equation}
with initial condition
\begin{equation}P_0h = h\end{equation}
It is a natural question to know if there is a heat kernel:
\begin{equation}P_th(x) = \int_Mp_t(x,y)h(y)dy\end{equation}
There are several ways to solve this problem:

-)The microlocal analysis ([12],[18],[19]) which uses as basic tool the Fourier transform and some regularity on the coefficients of $L$. In the case of a partial differential operator on $\mathbb{R}^d$, this means that $L = \sum a_{(\alpha)}(x){\partial^{(\alpha)}\over \partial x^{(\alpha)}}$ where $(\alpha)$ is a multiindex and $x \rightarrow a_{(\alpha}(x)$ is smooth.

-)The harmonic analysis, which uses as basic tools functional inequalities and does not need any regularity on the coefficients of $L$ ([3], [13], [51]).

-)The Malliavin Calculus([20], [44], [49]), which works for Markov semi-groups: $P_tf \geq 0$ if $f\geq 0$. The Malliavin Calculus requires moreover that the semi-group is  represented by a stochastic differential equation.

More precisely, the Malliavin Calculus needs a probabilistic representation of the semi-group $P_t$ by using the theory of stochastic differential equations where a flat Brownian motion or a Poisson process play a fundamental role.

Let us recall the main idea of the Malliavin Calculus in the case of the flat Brownian motion. Let us consider the Hilbert space $\mathbb{H}$ of finite energy maps starting from 0 from $[0,1]$ into $\mathbb{R}^m$ $t \rightarrow r_t = (r^i_t)$ endowed with the Hilbert norm
\begin{equation}\Vert r \Vert^2 = \sum_{i=1}^m \int_0^1\vert d/dtr_t^i\vert^2dt\end{equation}
We consider the formal Gaussian measure on $\mathbb{H}$ (written in the heuristic way of Feynman path integral)
\begin{equation}d\mu(r) = 1/Z\exp|-\Vert r\Vert^2/2]dD(r)\end{equation}
where $dD(r)$ is the formal Lebesgue measure on $\mathbb{H}$. Haar measure satisfying all the axioms of measure theory on a group exists if and only the group is locally compact. (We refer to [2] and [30] to defined Haar measure in infinite dimensin in a generalized way). This explains that we need to construct this measure on a bigger space, the space of continuous fonction $C([0,1],\mathbb{R}^m)$$t \rightarrow B_t$  issued from $0$  from $[0,1]$ into $\mathbb{R}^m$. There are a lot of Gaussian measures on $C([0,1],\mathbb{R}^m)$ ([48]]) but  the law of the Brownian motion is related to the heat equation on $\mathbb{R}^m$
\begin{equation}{\partial \over \partial t}P_tf(x) = 1/2\sum_{i=1}^m{\partial^2\over \partial x_i^2}P_tf(x)\end{equation}
We have namely
\begin{equation}P_th(x) = E[h(B_t+x)]\end{equation}
if $f$ is a bounded continuous function on $\mathbb{R}^m$. In such a case we have a semigroup operating on continuous function on $\mathbb{R}^m$.

We consider m smooth vector fields on $\mathbb{R}^d$ with bounded derivatives at each order. Vector fields here are considered as first order partial differential operators. We consider the operator
\begin{equation} L = 1/2\sum_{i=1}^mX_i^2\end{equation}
We introduce the Stratonovitch differential equation ([20], [49]) starting from $x$ (Vector fields here are considered as vectors whichs depends smoothly of $x$).:
\begin{equation}dx_t(x) = \sum_{i=1}^mX_i(x_t(x))dB_t^i\end{equation}
This is (and not the Ito equation) the correct equation associated to
\begin{equation}dx_t(r)(x) = \sum_{i=1}^mX_i(x_t(h)(x))dr^i_t\end{equation}
for $r \in \mathbb{H}$ endowed with the formal Gaussian measure $d\mu(r)$.

By Ito Calculus ([20],[49]), we can show that the semigroup $P_t$ generated by $L = 1/2\sum_{i=1}^mX_i^2$ is related to the diffusion $x_t(x)$ by the formula
\begin{equation}P_t(h)(x) = E[h(x_t(x))]\end{equation}
if $h$ is a continuous function on $\mathbb{R}^d$ (In such a case, the semigroup acts on continuous bounded functions on $\mathbb{R}^d$).

Malliavin idea is the following([44]): he differentiates in a generalized sense the It\^o  map $B_. \rightarrow x_t(x)$ . If this Ito map is a submersion in a generalized sense (The inverse of the Malliavin matrix belongs to all the $L^p$), the law of $x_t(x)$ has a smooth density and therefore the semigroup has an heat kernel. Malliavin for that uses a  heavy apparatus  of differential  operations on the Wiener space. Let us recall that there are several pioneering works of the Malliavin Calculus ([1], [6], [16]) motivated by mathematical physics, but only Malliavin Calculus is adapted to the study of stochastic differential equations and fit very well to the study of all measures of stochastic analysis.

Bismut ([7]) avoids to use this  heavy apparatus of differential operations on the Wiener space, by using a suitable Girsanov transformation and a system of convenient stochastic differential equations in cascade associated to the original stochastic differential equation. This allows to Bismut's way to get in a simpler way the Malliavin integration by parts for diffusions: if $(\alpha)$ is a multiindex, if $t>0$,
\begin{equation}E[h^{(\alpha)}x_t(x))] =  E[h(x_t(x))Q_t^{(\alpha)}]\end{equation}
where $Q_t^{(\alpha)}$ is a polynomial in the extra compoents of the system of stochastic differential equations in cascade and in the inverse of the Malliavin matrix.

The fact that only stochastic differential equations in cascade (therefore a system of semi-groups in cascade) appear in Bismut's approach of the Malliavin Calculus allows us to interpret Bismut's way of the Malliavin Calculus in the theory of semigroup by expulsating the probabilist language in [31]. We refer to [32], [33] for reviews with some applications.

[31] uses an elementary integration by parts, which has to be optimized. The main remark is that we can adapt this elementary integration by parts for non-markovian semi-groups.. It is possible to adapts Bismut's way of the Malliavin Calculus for non-markovian semi-groups.

It is divided in two steps:

-)An algebra on the semi-group. Only existences on the semigroup are required.

-)Estimates on the enlarged semigroup, which are necessary because polynomial function appear in the Malliavin integration by parts which are not bounded, which are performed in the non-markovian case by the Davies gauge transform (In the Markovian case, they were done by an adaptation in semi-group on the classic Burkholder-Davies-Gundy inequalities of stochastic analysis).

Moreover, Bismut in his seminal work ([9]) has done an intrinsic integration by part formula for the Brownian motion on a manifold, which overcame the problem that in the standard Malliavin Calculus there are a lot of stochastic  differential equations which represent the {\bf{same}} semigroup. In part IV we perform an intrinsic Malliavin Calculus  associated to a wide class of pseudodifferential elliptic operator, by performing a variation of the original pseudodifferential operator by a fractional power of it {\bf{intrisically}} associated to the original operator. We do the relation between the Malliavin Calculus of Bismut type and the general theory of elliptic pseudodifferential operators.

Bismut in his seminal work [9] pointed out the relation between the Malliavin Calculus and the large deviation theory for the study of short time asymptotics of the heat-kernel associated to diffusion semi-groups. We refer to the reviews [26], [29], [53], the book [5] and the seminal work [47] for probabilist methods in short time asymptotics of semi-groups.

Let us recall quickly the main goal of large deviation theory, here of Wentzel-Freidlin type [4], [52] and [54]. We introduce a small parameter and consider the stochastic differential equation with a small parameter starting from $x$:
\begin{equation}dx_t^\epsilon(x) = \epsilon\sum_{i=1}^mX_i(x_t^\epsilon)(x)dB_t^i\end{equation}
Wentzel-Freidlin theory  allows to get estimates of the type, when $\epsilon \rightarrow 0$
\begin{equation}\lim 2\epsilon^2Log[P[x_.^\epsilon(x) \in 0] = -\inf_{x_.(h)(x) \in O} \Vert r \Vert^2\end{equation}
if $O$ is an open subset of $C([0,1], \mathbb{R}^d)$ equipped with the uniform norm. We don't give details of the lot of technicalities in this estimate.

It is possible to adapt ([35], [37], [38], [39], [40]) Wentzel-Freidlin estimates to the case of non-markovian semi groups with the normalisation of W.K.B. analysis  of Maslov school ([45]) (See [17], [27] for seminal works on W.K.B. analysis). The main remark is that we can get only upper-bounds, because the semi-group does not preserves the positivity in this case. The second remark is that these estimates are valid only for the semi-group, because in this case path space functional integrals are not defined (See [36] for a review and the work [11], [25], [46]). The normalizations are those classical of semi-classical analysis but the type of estimates is different. They work for the heat equation and not for the Schroedinger equation.

This allows to fullfill in this non-markovian context the beautifull request of Bismut's book [5] and to do the marriage between the Malliavin Calculus and Wentzel-Freidlin estimates. The main difference is that we have to consider the absolute value of the heat-kernel because in such a case the semi-group does not preserve the positivity such that we get only upper-bound in the studied Varadhan type estimates (Wentzel-Freidlin estimates are still valid for the heat-kernel).

This work is a review paper of several of our works. The main novelty is part IV, which is new.
\section{The case of a formal stochastic differential equation}

 Let us consider an elliptic differential operator of order $l$ on a compact manifold $M$ of dimension $d$. If we perturb it by a strictly lower order operator $L_p$, it results by 
the theory of pseudo-differential operator (which is given by the role of the principal symbol of an elliptic operator)  that the qualitative behaviour (hypoellipticity..) is the same than the qualitative behavior of $L+L_p$. See [12],[18],[19] for various textbooks in analysis about this problematic.

Recently, we have introduced an elliptic operator of order $2k$ $L_0= \sum f_i^{2k}$ where $f_i$ is an orthonormal basis of the Lie algebra of a compact Lie group $G$ of dimension $m$ with generic element $g$.  $f_i$  are considered 
as right invariant vector fields.We have established the Malliavin Calculus of Bismut type for $L$.  We consider a polynomial $Q$ of degree strictly smaller than $2k$  in the vector fields$ f_i$ with constant components. We consider the total operator
\begin{equation} L = L_0+Q\end{equation}
The goal of this part, by using a small interpretation of [41] and [42]  is to adapt in this present situation the strategy of [41] for diffusions.. ( [41] [42] have used the machinery of the Malliavin Calculus [7] translated by ourself in semi-group theory for diffusions in [31]) Malliavin matrix plays here a fundamental role in the optimization of the integration by parts in order to arrive to full Malliavin integration by parts. All formulas are {\bf{formally the same}} if we add or not add the perturbation of the main operator..

We consider the elliptic operator on $G \times \mathbb{R}$
\begin{equation}Q + \sum_i^df_i^{2k} + \sum r_{i,t}f_i{\partial\over \partial u} + {\partial^{2k}\over \partial u^{2k}} = \tilde{L}_t^r\end{equation}
It generates by elliptic theory a semi-group   on $C_b(G\times \mathbb{R})$, the space of bounded continuous function on $G\times \mathbb{R}$ endowed with the uniform norm..
\begin{theorem} (Elementary integration by parts formula).We have if $h$ is smooth with compact support
\begin{equation}\int_0^tP_{t-s}\sum h_{s,i}e_iP_s[h]ds = \tilde{P}_t^h[uh](.,0)\end{equation}\end{theorem}
{\bf{Proof}}: It is the same proof than the proof of Theorem 3 of [42]. $\diamondsuit$

Let $ V= G \times M_d$. $M_d$ is the space of symmetric matrices on 
$Lie G$. $(x,v) \in V$. $v$ is called the Malliavin matrix. We consider
\begin{equation}\hat{X}_0 = (0, \sum <g^{-1}f_i,.>^2)\end{equation}

We consider the Malliavin generator (We skipp the problems of signs)
\begin{equation}\hat{L} = Q + \sum f_i^{2k}-\hat{X}_0\end{equation}
\begin{theorem}$\hat{L}$ spanns a semi-group.  $\hat{P}_t$ called the Malliavin semi-group on $C_b(M)$.\end{theorem}
{\bf{proof}}  It is the same proof of theorem 4 of [42] since $Q$ is a polynomial with constant compoents in the $f_i$ and $L$ generates a $C_b(G)$ semi-group. The proof leads to some difficulties because the Malliavin operator is not the perturbation of an elliptic operator and uses the Volterra expansion.
$\diamondsuit$

The Malliavin semi-group will allow to us to optimize the elementary integration by parts of theorem 2.
We have the main theorem of this paper:
\begin{theorem}(Malliavin) If the Malliavin condition holds
\begin{equation}\vert\hat{P}_t][v^{-p}](g,0) < \infty\end{equation}
for all integer positive integer $p$, $P_t$ has an heat-kernel.\end{theorem}
{\bf{Proof}}: It is the same proof as in the beginning of the proof of theorem 6 of [42]. Under Malliavin assumption, we can optimize the elementary integration by part of Theorem 2, in order to get, according the framework of the Malliavin Calculus, the inequality for any smooth function $h$ on $G$
\begin{equation}\vert P_t[<dh,f_i>]\vert \leq C \Vert h \Vert_\infty\end{equation}
$\diamondsuit$

{\bf{Remark}}: Let us explain quickly the philosophy of this theorem, when there is no perturbation term.
We consider  a set of path in $\mathbb{R}^d$ denoted  $r^i_t$ which represent the semi-group asoociated to $\sum_i^k {\partial^{2k}\over \partial u_i^{2k}}$. We don't enter in the problem of signs. We consider the formal stochastic differential equation
\begin{equation}dx_t(r)(e) = \sum_i^df_idr_t^i\end{equation}
issued from  $e$. Formally, this represent the semi-group $P_t$ without the perturbation term
\begin{equation}P_t[h](e) =  "E"[f(x_t(e)]\end{equation}
Malliavin assumption express in some sense that the "Ito" map $r_.^. \rightarrow x_t(e)$ is a submersion.

By this inequality, we deduce according the framework of the Malliavin Calculus
that
\begin{equation}P_t[h](e) = \int_Gh(g)p_t(e,g)dg\end{equation} for a non strictly positive heat-kernel $p_t$ ($dg$) denotes the normalized  Haar mesure on $G$).
if the Malliavin assumption is satisfied.

\begin{theorem}Under the previous elliptic assumptions, 
\begin{equation}\vert \hat{P}_t \vert [\vert v^{-p}\vert]](g_0,0) < \infty\end{equation} if $t>0$\end{theorem}

{\bf{Proof}} It is the same proof than the proof of theorem 8 of [42].It is based upon the initial strategy to invert the Malliavin matrix in stochastic analysis by slicing the time interval in small time intervals . Only the main part of the generator plays the main role in this strategy because we are in an elliptic case.$\diamondsuit$

We can iterate the integration by parts formulas, by introducing a system of semi-groups in cascade. We deduce the theorem
\begin{theorem}If $t>0$ the semi group $P_t$ has a smooth heat kernel
\begin{equation}P_t([h](g) = \int_Gp_t(g,g')dg'\end{equation}\end{theorem}
The main remark is that the heat kernel can change of sign. This theorem is classical in analysis [51] but it enters in our general strategy to implement stochastic tools in the general theory of linear semi-groups.

In order to simplify the computation, we have used the symmetry of the group. In the next part, we will use fully the symmetry of the group to simplify the computations.
\section{The full use of the symmetry of the group}
Let us recall what is a pseudodifferential operator on $\mathbb{R}^d$ ([12], [17], [18]). Let be a smooth  function function  from  $\mathbb{R}^d \times \mathbb{R}^d$ into $\mathbb{R}$ $a(x, \xi)$. We suppose that
\begin{equation}\sup_{x \in \mathbb{R}^d}\vert  D_x^rD_\xi^la(x,\xi)\vert \leq C \vert \xi\vert^{m-l}+C\end{equation}
We suppose that
\begin{equation}\inf_{x \in \mathbb{R}^d}\vert a(x, \xi)\vert \geq C \vert \xi \vert^{m'}\end{equation}
 for $\vert \xi \vert > C$
 for a suitable $m'>0$.
Let $\hat{h}$ the fourier transform of the continuous function $h$. We consider the operator $L$ defines on smooth function $h$ by
:\begin{equation}\hat{Lh}(x) = \int_{\mathbb{R}^d}a(x,\xi)\hat{h}(\xi)d\xi\end{equation}
$L$ is said to be a pesudodifferential operator elliptic of order larger than $m'$ with symbol $a$. This property is invariant if we do a diffeomorphism on $\mathbb{R}^d$ with bounded derivatives at each order. This remark allows to define by using charts  a pseudodifferential operator elliptic of order larger than $m'$ on a compact manifold $M$.

 Let $f^i$ be a basis of $T_eG$. We can consider as rightinvariant vector fields. This means that if we consider the action $R_{g_0}$ $h \rightarrow (g\rightarrow h(gg_0))$
on smooth function $h$ on $G$, we have
\begin{equation}R_{g_0}(f^ih) = f^i(R_{g_0}h)\end{equation}.

We consider a rightinvariant elliptic pseudodifferential positive elliptic operator $L$ of order larger than $2k$ on $G$ . It generates by elliptic theory a semi group $P_t$ on $L^2(dg)$ and even on $C_b(G)$ the space of continuous functions on $G$ endowed with the uniform norm. 

\begin{theorem}If $t>0$,
\begin{equation}P_th(g_0) = \int_Gp_t(g_0,g)h(g)dg\end{equation}
where $g \rightarrow p_t(g_0,g)$ is smooth if $h$ is continuous.   \end{theorem}

This theorem is classical in analysis , but it enters in our general program to implement stochastic analysis tool in the theory of Non-Markovian semi-group. See the review [36] for that. See [41], [42] for another presentation where the Malliavin Matrix plays a key role. Here we don't use the Malliavin matrix. See [43] fot the case of rightinvariant differential operators.
 The proof is divided in two steps.
\subsection{Algebraic scheme of the proof: Malliavin integration by parts}
We consider the family of operators on $C^\infty(G\times \mathbb{R}^n)$:
\begin{equation}\tilde{L}_t^n = L + \sum_{i=1}^n f^{j_i}{\partial \over \partial u_i}\alpha_t^i + \sum_{i=1}^n{\partial^{2k}\over \partial u_i^{2k}}\end{equation}
$\alpha_t^i$ are smooth function from $\mathbb{R}^+$ into $\mathbb{R}$.
By elliptic theory, $\tilde{L}_t^n$ generates  a semi-group $\tilde{P}_t^n$ on $C_b(G\times \mathbb{R}^n)$. This semi-group is time inhomegeneous.
\begin{equation}\tilde{P}_t^{n+1}[h(g)h^n(u)v](.,.,0) =\\ \int_0^t\tilde{P}_{t,s}^n[f^{j_{n+1}}\alpha_s^{n+1}\tilde{P}^n_s[h(g)h^n(u)](.,.)\end{equation}
Moreover
\begin{equation}\tilde{P}_t^{n+1}[uh(.)h^n(.)](.,.,u_{n+1}) = \\ \tilde{P}_t^{n+1}[uh(.)h^n(.)](.,.,0)+ \tilde{P}_t^{n}[h(.)h^n(.)](.,.)u_{n+1} \end{equation}
$h$ is a function of $g$, $h^n$ a function of $u_1,..., u_n$. This comes from the fact that  ${\partial\over \partial u_{n+1}}$ commute with the considered operator.

Therefore the two sides of (37) satisfy the same parabolic equation with second-member.
We deduce that
\begin{equation}\tilde{P}_t^{n+1}[u_{n+1}\prod_{j=1}^nu_jh(.)](.,.,0) = \\ \int_0^tds\tilde{P}_{t,s}^n[f^{j_{n+1}}\alpha_s^{n+1}\tilde{P}_s^n[h\prod_{j=1}^nu_j]](.,.)\end{equation}
This is an integration by parts formula. We would like to present this formula in a more appropriate way for our object.

We consider the operator 
\begin{equation} \overline{L}^n  =  L + \sum_{j=1}^n {\partial^{2k}\over \partial u_j^{2k}}\end{equation}
It generates a semi-group $\overline{P}^n_t$. In the sequel we will skip the problem of sign coming if $k$ is even or not.

We introduce a suitable generator
\begin{equation}\tilde{R}_t^{n+1} = \overline{L}^n+ F_s\end{equation}
by taking care of the relation $[f^i,f^j] = \sum_k \lambda_k^{i,j}f^k$.
It is an operator of the type studied. It generates therefore a time inhomogeneous semi-group $\tilde{Q}_t^n$.  
Therefore the integration by parts formula (39) can be written in the more suitable way
\begin{multline}\tilde{P}_t^{n+1}[u_{n+1}\prod_{j=1}^nu_jh(.)](.,.,0) =  \int_0^t\alpha_s^{n+1}ds\tilde{P}_t^n[f^{j_{n+1}}h\prod_{i=1}^nu_i](.,.) +\\ \int_0^t\alpha_s^{n+1} ds \tilde{P}_{t,s}^n \tilde{Q}_s^n[h\prod_{i=1}^nu_i](.,.)\end{multline}
We do the following recursion hypothesis on $l$:

{\bf{Hypothesis (l)}} There exists a positive real $r_l$ such that if $(\alpha)$ is a multiindex of length smaller than l
\begin{equation}\vert \tilde{P}_t^n[f^{(\alpha)}h\prod_{i=n}^nu_i](g,v_.)\vert \leq Ct^{-r_l} \\ \Vert h\Vert_\infty (1 + \prod_{i=n}^n\vert v_i\vert )\end{equation}
where $\Vert . \Vert_\infty$ is the uniform norm of $h$.

It is true for $l=1$ by (39) and the estimates which follow.

If it is true for $l$, it is still true for $l+1$, by using (42) for $f^{(\alpha)}h$ and taking $\alpha_s^{n+1} = s^{r_l}$

By choosing suitable $\alpha_t^j$, we have according the framework of the Malliavin Calculus for any multiindex $(\alpha)$
\begin{equation}\vert P_t[f^{(\alpha)}h](g_0)\vert \leq C_{(\alpha)}\Vert h\Vert_\infty\end{equation}
in order to conclude.

 \subsection{Estimates: the Davies gauge transform}

We do as in [43] (26). The problem is that in $\tilde{P}_t^n[h\prod_{j=1}^nu_j](.,.)$ the test function $u_j$ are not bounded and that $\tilde{P}_t^n$ acts only on $C_b(G\times \mathbb{R}^n)$. We do as in [3] the Davies gauge transform $\prod_I^n g(u_i)$ where
\begin{equation}g(u) = (\vert u \vert )\end{equation}
if $u$ is big and $g$ is smooth .

This  gauge transform acts on the original operator by the simple  formula $(\prod_{i=1}^ng(u_i))^{-1}\tilde{L}_1^n((\prod_{i=1}^ng(u_i).)$. On the semi group it acts as
\begin{equation}(\prod_{i=1}^ng(.))^{-1}\tilde{P}_t^n[(\prod_{i=1}^ng(u_i)h(.)h^n(.)](.,.)\end{equation}
But 
\begin{equation}(g(u_i))^{-1}{\partial\over \partial u_i}(g(u_i).) = {\partial \over \partial u_i} + C(u_i)\end{equation}
where the potential $C(u_i)$ is smooth with bounded derivatives at each order. Therefore the transformed semi-group act on $C_b(G\times{R}^n)$. 

{\bf{Remark}}We can consider as particular case ([43])
Let $G$ be a compact connected Lie group, with generic element $g$ endowed with its binvariant Riemannian structure and with its normalized Haar measure $dg$. $e$ is the unit element of $G$.

Let $f^i$ be a basis of $T_eG$. We can consider as rightinvariant vector fields. This means that if we consider the action $R_{g_0}$ $h \rightarrow (g\rightarrow h(gg_0))$
on smooth function $h$ on $G$, we haver
\begin{equation}R_{g_0}(f^ih) = f^i(R_{g_0}h)\end{equation}.

Let be $\xi^{(\alpha)} = \xi^{\alpha_1}.... \xi^{\alpha^{\vert \alpha\vert}}$ and let be $f^{(\alpha)} = f^{\alpha_1}..f^{\alpha_{\vert \alpha\vert}}$. 
$(\alpha)$ is a multi-index of length $\vert \alpha\vert$.

We consider a matrix $a_{\alpha, \beta}$ for multindices of length $k$, which is supposed symmetric strictly positive.

We consider the operator
\begin{equation}L = \sum_{(\alpha), (\beta)}f^{(\alpha)}a_{(\alpha), ( \beta)}f^{(\beta)}\end{equation}
According [51], $(-1)^k L$  is a  positive symmetric densely elliptic defined operator on $L^2(G)$, which generates by elliptic theory a semi-group acting on $C_b(G)$, the space of continuous function on $G$. In such a case, we have an heat-kernel associated to the semi-group (See [43]). The case of a rightinvariant differential operator has exactly the same proof than the case of theorem 6, where the details will be presented elsewhere.
\section{The case of an intrinsic variation}
Let $L$ be a strictly positive self-adjoint operator on a compact manifold $M$. We suppose that $L$ is a pseudo-differential elliptic operator of order $l \geq 2k$ for an integer $k \geq 1$. It generates a contraction semi-group on $L^2(M)$ and by ellipticity a semi-group on $C_b(M)$.
\begin{theorem}There is an heat-kernel $p_t(x,y)$ associated to $P_t$. If $t>0$
\begin{equation}P_t(f)(x) = \int_Mp_t(x,y)f(y)dy\end{equation}
where $y \rightarrow p_t(x,y)$ is smooth.
\end{theorem}
The proof is divided in two steps:

\subsection {Algebraic scheme of the proof: Malliavin integration by parts}Let $\alpha$ belonging to $]0,1[$. The fractional power [50] $L^\alpha$ is still a strictly positive pseudodifferential operator elliptic of order $\alpha l$, which commutes with $L$. We skipp up later the problem if $k$ is even or not. We consider the operator on $C^\infty(M\times \mathbb{R}^n)$

\begin{equation}\tilde{L}^n_s = L + s^rL^\alpha \sum_i^n{\partial\over \partial u_i}+  \sum_{i=1}^n {\partial^{2k}\over \partial u_i^{2k}}\end{equation}
It is an elliptic operator of order $2k$ on $M\times \mathbb{R}^n$. The main part
\begin{equation}\overline{L}^n =  L +   \sum_{i=1}^n {\partial^{2k}\over \partial u_i^{2k}}\end{equation}
is positive and is essentially self-adjoint. Therefore the main part generates a semi-group on $C_b(M\times \mathbb{R}^n)$. This remains true for
 $\tilde{L}^n$ because $\tilde{L}^n$ is a perturbation of $\overline{L}^n$ by a strictly lower operator. We call this semi-group $\tilde{P}_t^n$.

The main remark is that $L^\alpha$ commutes with $\tilde{L}^n$ such that
\begin{equation} L^\alpha \tilde{P}_t^n = \tilde{P}_t^n  L^\alpha\end{equation}
According the beginning of the previous part, we get the elementary integration by part
\begin{multline}\tilde{P}_t^{n+1}[f\prod_{i=1}^nu_i u ](x, v_i,0)  = \int_0^tP_{t-s}^n[s^rL^{\alpha}\tilde{P}_s^n[f\prod_{i=1}^nu_i]](x,v_i) = \\ \tilde{P}_t^n[L^\alpha f \prod_{i=1}^nu_i](x,u_i)\int_0^ts^rds\end{multline}
Suppose by induction on $l$ that
\begin{equation}\vert \tilde{P}_t^n[(L^\alpha)^lf\prod_{i=1}^nu_i](x,v_i)\vert \leq  C t^{-r(l)}\Vert f \Vert_\infty(1 + \prod_{i=1}^n \vert v_i \vert ) \end{equation}
By applying the elementary integration by parts (54) to $(L^\alpha)^l)f$, and choosing $r = r(l)$, we deduce our result.
Therefore we have the inequality
\begin{equation}\vert P_t[(L^\alpha)^lf](x)\vert \leq C t^{-r(l)}\Vert f \Vert_\infty\end{equation}
The result follows from the fact that $L^\alpha$ is an elliptic operator.

\subsection{Estimates: the Davies gauge transform}
 We do as in [43] (26). The problem is that in $\tilde{P}_t^n[h\prod_{j=1}^nu_j](.,.)$ the test function $u_j$ are not bounded and that $\tilde{P}_t^n$ acts only on $C_b(G\times \mathbb{R}^n)$. We do as in [35] the Davies gauge transform $\prod_I^n g(u_i)$ where
\begin{equation}g(u) = (\vert u \vert )\end{equation}
if $u$ is big and $g$ is smooth strictly positive .

This  gauge transform acts on the original operator by the simple  formula $(\prod_{i=1}^ng(u_i))^{-1}\tilde{L}_1^n((\prod_{i=1}^ng(u_i).)$. On the semi group it acts as
\begin{equation}(\prod_{i=1}^ng(.))^{-1}\tilde{P}_t^n[(\prod_{i=1}^ng(u_i)h(.)h^n(.)](.,.)\end{equation}
But 
\begin{equation}(g(u_i))^{-1}{\partial\over \partial u_i}(g(u_i).) = {\partial \over \partial u_i} + C(u_i)\end{equation}
where the potential $C(u_i)$ is smooth with bounded derivatives at each order. Therefore the transformed semi-group act on $C_b(G\times{R}^n)$. It remains to choose
\begin{equation}h^n(u_.) =\prod_{j=1}^n{u_j\over g(u_j)}\end{equation}
in order to conclude. We deduce the bound:

\begin{equation}\vert \tilde{P}_t^n\vert[h\prod_{j=1}^n\vert u_j\vert](.; v_.) \leq C(\Vert h\Vert_\infty (1 + \prod_{i=n}^n\vert v_i \vert)\end{equation}
where $\vert \tilde{P}_t^n\vert$ is the absolute value of the semi-group $\tilde{P}_t^n$.

{\bf{Remark:}}We could show that $(x,y) \rightarrow p_t(x,y)$ is smooth if $t>0$ by the same argument.

{\bf{Remark:}}We can replace the hypothesis $L$ strictly positive by the hypothesis $L$ positive by replacing $L^\alpha$ by $(L+CI_d)^\alpha$ where $C>0$.

\section{Wentzel-Freidlin estimates for the semi-group only}
We consider a differential operator of order $2k$ on the compact manifold $M$  which is supposed elliptic of order $2k$ and strictly positive. We suppose we can write it
as
\begin{equation} L= \sum_{j=0}^{2k}\sum_{i=0}^{r(j)}(X_{i,j})^{j}\end{equation}
where $X_{i,j}$ are smooth vector fields on $M$. The ellipticity assumption states that
\begin{equation}\sum_{i=0}^{r(2k)}<X_{i,2k},\xi>^{2k} = H(x,\xi)  \geq C \vert \xi\vert^{2k}\end{equation}
To the Hamiltonian $H$, we introduce the Lagrangian
\begin{equation}L(x,p) = \sup_\xi(<p,\xi>-H(x,\xi))\end{equation}
We get the estimate
\begin{equation}-C + C\vert p \vert^{{2k\over 2k-1}}\leq L(x,p) \leq C +\vert p \vert^{{2k\over 2k-1}}\end{equation}
for some strictly positive constants $C$.

If $\phi$ is a continuous piecewise differentiable path on $M$, we put:

\begin{equation} S(\phi) = \int_0^1L(\phi(t), d/dt\phi(t))dt \end{equation}
and we put
\begin{equation} l(x,y) =  \inf_{\phi(0) = x, \phi(1) = y}S(\phi)\end{equation}
By Ascoli theorem, $(x,y) \rightarrow l(x,y)$ is a continuous function on $M\times M$.
\begin{theorem}(Wentzel-Freidlin)If $O$ is an open ball of $M$, we have when $t \rightarrow 0$
\begin{equation}\overline{\lim}t^{{1\over 2k-1}}\log\vert P_t\vert(1_O](x) \leq-\inf_{y \in O}l(x,y)\end{equation} \end{theorem}
{\it{Proof:}}We put $\epsilon = t^{{1\over 2k-1}}$. According the normalisation of Maslov school [37], we consider the semi-group $P_s^\epsilon$ associated to $L_\epsilon= \epsilon^{2k-1}L$. Moreover
\begin{equation}P_t = P_1^t\end{equation}
where 
$P_s^t$ is associated to $tL$ ([10])
The result will arise if we show when $\epsilon \rightarrow 0$

\begin{equation}\overline{\lim} \epsilon\log\vert P_1^\epsilon\vert(1_O](x) \leq -\inf_{y \in O}l(x,y)\end{equation}
The main ingredient is:
\begin{lemma}For all $\delta>0$, all $C$, there exists $s_\delta$ such that if $s< s_\delta$
\begin{equation}\vert P_s^\epsilon\vert[1_{B(x,\delta)^c}](x) \leq \exp[-C/\epsilon]\end{equation}
where $B(x, \delta)$ is the ball of radius $\delta$ and center $x$\end{lemma}
{\it{Proof of the lemma}} We imbedd $M$ in a linear space. We consider the semi-group
\begin{equation}Q_s^\epsilon = \exp[-<x,\xi>/\epsilon]P_s^\epsilon[ \exp[<x',\xi>/\epsilon](x').](x) \end{equation}
Its generator is 
\begin{equation}\overline{L}_\epsilon + H(x, \xi)/\epsilon\end{equation}
\begin{equation}\overline{L}_\epsilon = L_\epsilon + R_\epsilon\end{equation}
In the perturbation term $R_\epsilon$, there are only differential operators of order $l$, $l\in]0,2k[$. When a differential operator of degree l appears, there is a power of at least $l-1$ of $\epsilon$ which appears and a power of $\xi$ at most $2k$ which appears.

Let us consider in a small neighborhhod of x the diffeomorphism
\begin{equation}\Psi_\epsilon: y \rightarrow x + {y-x\over \epsilon^{{2k-1\over 2k}}}\end{equation}
Outside a big neighborhood of $x$, $\Psi_\epsilon$ is the identity.

We consider the measure $\mu_\epsilon$
\begin{equation}f \rightarrow P_1^\epsilon[F(\Psi_\epsilon(x))](x)\end{equation}
Under the transformation $\Psi_\epsilon$, the vector fields $\epsilon^{{2k-1\over 2k}} X_{i,j}$ are transformed in the vector field $X_{i,j}(x+  \epsilon^{{2k-1\over 2k}}(y-x))$. Therefore we can apply the machinery of the previous part in order to show that the measure $\mu_\epsilon$ has a bounded density $q_\epsilon(x,.)$ when $\epsilon \rightarrow 0$.

Let $R$ be a differential operator of order $l$.
 We have
\begin{equation}\int_Mg(x)RP_1^\epsilon[h](x)dx = \int_{M\times M}g(x)h(y)R_xp_1^{\epsilon}(x,y)dxdy\end{equation}
By symmetry
\begin{equation}p_1^\epsilon(x,y) = p_1^\epsilon(y,x)\end{equation}
 Then

\begin{equation}\int_Mg(x)RP_1^\epsilon[h](x)dx = \int_M h(y)P_1^\epsilon[Rg](y)dy\end{equation}
By the previous remark
\begin{equation} \vert P_1^\epsilon[Rh](y)\vert \leq {C\over  \epsilon^{l{2k-1\over 2k}}}\Vert h\Vert_\infty\end{equation}
Therefore
\begin{equation} \vert \int_Mg(x)RP_1^\epsilon[h](x)dx  \vert \leq
{C\over  \epsilon^{l{2k-1\over 2k}}}\Vert g\Vert_\infty \Vert h \Vert_\infty\end{equation}
We deduce that
\begin{equation} \vert RP_1^\epsilon[h](x)\vert   \leq {C\over  \epsilon^{l{2k-1\over 2k}}}\Vert h\Vert_\infty\end{equation}
We deduce a bound of $R_\epsilon P_s^\epsilon$
\begin{equation}\vert R_\epsilon P_s^\epsilon h(x)\vert \leq {\vert \xi \vert^{2k-1}\over s^{{l\over 2k}}}\epsilon^{-1+1/k}\Vert h\Vert_\infty\end{equation}
We apply Volterra  expansion to $Q_s ^\epsilon$. We get
\begin{equation}\vert Q_s^\epsilon f\vert \leq \vert P_s^\epsilon f\vert  + \sum_{i=1}^\infty\vert \int_{\Delta_l(s)}I_{s_1,..,s_l}ds_1...ds_l\vert \end{equation}
where $\Delta_l(s)$ is the simplex $0<s_1<..<s_l<s$ and
\begin{equation}I_{s_1,..,s_l} = P^\epsilon_{s_1} (R_\epsilon + H/\epsilon)...P_{s_{l}-s_{l-1}}^\epsilon (R_\epsilon + H/\epsilon)P^\epsilon_{s-s_{l-1}} h\end{equation}
We deduce a bound of $\vert \int_{\Delta_l(s)}I_{s_1,..,s_l}ds_1...ds_l\vert $ by
\begin{equation}{\vert \xi \vert^{2lk}\over \epsilon^l}\int_{\Delta_l(s)}\prod (s_{i+1}-s_i)^{-{2k-1\over 2k}}ds_1..ds_l = {\vert \xi \vert^{2lk}\over \epsilon^l} I_l(s)\end{equation}
We suppose  by induction that
\begin{equation}I_{l}(s) = \alpha_ ls^{l(1+\beta_k)}\end{equation}
where $\beta_k \in ]-1,0[$. It is still true by the recurtion formula
\begin{equation}I_{l+1}(s) = \int_0^sI_l(u)  (s-u)^{-{2k-1\over 2k}}du\end{equation}
We deduce the bound 
\begin{equation}
\alpha_l \leq {C^l\over l!}\end{equation}
Therefore
\begin{equation} \vert Q_s^\epsilon h(x)\vert  \leq \exp[Cs\vert \xi\vert^{2k}/\epsilon] \Vert h \Vert_\infty\end{equation}
It remains to remark that we have the bound
\begin{equation}\vert P_s^\xi\vert[1_{B(x,\delta)^c}](x) \leq \exp[-{C\delta \vert \xi \vert \over \epsilon}+Cs\vert \xi\vert^{2k}/\epsilon] \end{equation}
and to extremize in $\vert \xi \vert$ to conclude.$\diamondsuit$

{\bf{End of the proof of theorem 9}}
We operate as in Freidlin-Wentzel book [54]  and as in [35],[38], and [39]. We slice the time interval $[0,,1]$ in a finite numbers of time intervals $[s_i, s_{i+1}]$ where we can apply the previous lemma. We deduce a positive measure on the set of polygonal paths, where we can repeat exactly the considerations of [35].

{\bf{Remark:}}This estimate is a semi-classical estimate with different type of estimates of W.K.B. estimates a la Maslov and with a different method. We consider in W.K.B. estimate a symbol of an operator $a(x, \xi)$ and we consider the generator $L_\epsilon$ associated with the normalized symbol (a la Maslov) $1/\epsilon a(x,\epsilon \xi)$. Let us suppose that $L_\epsilon$ generates a semi-group $P_t^\epsilon$. The object of WKB method is to get {\bf{precise}} estimates of the semi-group $P_1^\epsilon$ when $\epsilon \rightarrow 0$. For that people look at a formal asymptotic expansion (we omit to write the initial conditions) of $P_1^\epsilon$
of the type
\begin{equation}\epsilon^{-r}\exp[-l(y)/\epsilon]\sum\epsilon^iC_i(y)\end{equation}
The function $l$ satisfy a highly non-linear equation (the Hamilton-Jacobi-Belman equation) and $c_i(y)$ satisfy formally a system of linear partial differential equation in cascade. The cost function in theorem $l(x,y)$ is solution of the highly non-linear Hamilton-Jacobi-Belman equation, which is difficult to solve. Instead of precise asymptotics, we are interested by logarithmic estimates which are totally different with a method totally different. On the other hand, generally semi-classical asymptotics considers the case of the Schroedinger equation instead of the heat semi-group.

On $\mathbb{R}^d$ we can speak without any difficulty of the symbol of an operator. Poisson processes, L\' evy processes and jump processes are more or less generated by pseudodifferential operators whose generator satisfy the maximum principle (See [10], [21], [22], [13], [24], [28]). We will present pseudodifferential operators with a type of compensation of stochastic analysis which do not satisfy the maximum principle. The end of this part is extracted from [35] and [40].
Let us consider the generator on $C_\infty(\mathbb{R}^d)$
\begin{equation}Lf(x) = (-)^{l+1}\int_{\mathbb{R}^d}(f(x+y)-f(x)-\sum_{i=1}^{2l}<y^{\otimes i}, D^if(x)>){h(x,y)\over \vert y \vert^{2l+1+\alpha}}dy \end{equation}
$\alpha \in ]-1,0[$ $h(x,y) = 0$ if $\vert y \vert > C$ and $h\geq 0$.
The measure ${h(x,y)\over \vert y \vert^{2l+1+\alpha}}dy$ is called the L\'evy measure.
\begin{theorem}If $h(x,0) = 1$, $L$ is an elliptic pseudo-differential generator.\end{theorem}
\begin{definition}If $h(x,y) = h(y)$, we will say that $L$ is a generalized L\'evy generator.\end{definition}

\begin{theorem}Suppose that $L$ is of L\'evy type and that $h(y) = h(-y)$. $L$ is positive symmetic, and therefore admits by ellipticity a selfadjoint extension on $L^2(\mathbb{R}^d)$, which generates a contraction semi group on $L^2(\mathbb{R}^d)$ which is still a semi-group on $C_b(\mathbb{R}^d)$.\end{theorem}
{\bf{Remark:}}The symbol $a(x,\xi)$ of the generator is given by
\begin{equation}(-)^{l+1}\int_{\mathbb{R}^d}(\exp[\sqrt{-1}<y,\xi>]-\sum_{i=1}^{2l}({(\sqrt{-1}<\xi,y>)^i)\over i!}{h(x,y)\over \vert y \vert^{2l+1+\alpha}}dy \end{equation}
The Hamiltonian associated is the symbol in real phase. Let us conider a generator of L\'evy type of the previous theorem: it is
\begin{equation}(-)^{l+1}\int_{\mathbb{R}^d}(\exp[<y,\xi>]-\sum_{i=1}^{l}({<\xi,y>^{2i}\over{2 i}!}){h(,y)\over \vert y \vert^{2l+1+\alpha}}dy \end{equation}
The Hamiltonian is a smooth convex function equals to 1 in 0.
Associated to it, we consider the Lagrangian:
\begin{equation}L(p) = \sup_\xi (<\xi,p>-H(\xi))\end{equation}
If $t \rightarrow \phi_t$ is a piecewise differentiable continuous curve in $\mathbb{R}^d$, we consider its action $\int_0^1dtL(\phi_t, d/dt \phi_t) = S(\phi_.)$
We introduce the control function
\begin{equation}l(x,y) = \inf_{\phi_0= x; \phi_1 = y}S(\phi)\end{equation}
Let us recall that $(x,y) \rightarrow l(x,y)$ is positive finite continuous. 

We consider the generator associated to $1/\epsilon a(\epsilon \xi)$. This correponds in the classical case of jump process where the compensation is only of one term to the case of a jump process with more and more jumps which are more and more small [54].
We consider the generator $L^\epsilon$ associated to $1/\epsilon a(\epsilon \xi)$. It generates a semi-group $P_t^\epsilon$. We get:
\begin{theorem} Wentzel-Freidlin ([35], [40]). When $\epsilon \rightarrow 0$, we get if $O$ is an open ball of $\mathbb{R}^d$ if $l+1$ is even:
\begin{equation} \overline{\lim}\epsilon \log\vert P_1^\epsilon\vert[1_O](x) \leq - \inf_{y \in O}l(x,y)\end{equation}\end{theorem}
{\bf{Remark}}: For this type of operator, Wentzel-Freidlin estimates are not related to short time asymtotics.

\section{Application: some Varadhan estimates}
This part follows closely [43]. Only the mechanism of the integration by part is different from [39]. For large deviation estimates with respect of W.K.B normalization at the manner of Maslov [45] for Non-Markovian operators, we refer to [38] for instance.

Let us consider the Hamiltonian function from $T^*(G)$ into $\mathbb{R}^+$
\begin{equation}H(g,\xi) = \sum_{\vert \alpha \vert = k, \vert \beta \vert = k}\\<f^{(\alpha_1)},\xi>..<f^{(\alpha_k)},\xi>a_{(\alpha), (\beta)}<f^{(\beta_1)}..\xi>..<f^{(\beta_k)},\xi>\end{equation}
$H(g,p)$ is positive convex in $p$. According the theory of  large deviation, we consider the associated Lagrangian
\begin{equation}L(g,\xi) = \sup_\xi(<\xi,p>-H(g,\xi))\end{equation}
If $t \rightarrow \phi_t$ is a curve in the group, we consider its action $\int_0^1dtL(\phi_t, d/dt \phi_t) = S(\phi_.)$
We introduce the control function
\begin{equation}l(g_0,g_1) = \inf_{\phi_0= g_0; \phi_1 = g_1}S(\phi)\end{equation}
Let us recall that $(g_0,g_1) \rightarrow l(g_0,g_1)$ is positive finite continuous. 

We have shown in the previous part that  if we consider a small parameter $\epsilon$ and if we consider the generator $\epsilon^{2k-1}L$ and the semi group $P_t^\epsilon$ associated and if $g_0$ and $g_1$  are not closed , we get for any small ball centered in $g_1$ uniformly:

\begin{equation}\overline{Lim}_{\epsilon \rightarrow 0}\epsilon Log\vert P_1^\epsilon \vert[1_O](g_0) \leq -\inf_{g_1 \in O}l(g_0,g_1)\end{equation}
where$\vert P_1^\epsilon \vert$ is the absolute value of the semi-group (See [38]). See for that the previous part

But $P_t = P_1^t$ where $P_s^t$ is the semi group associated to $tL$ (See [15]). We put $\epsilon = t^{1/2k-1}$ such that 
\begin{equation}\overline{Lim}_{t \rightarrow 0}t^{1/2k-1} Log\vert P_t \vert[1_O](g_0) \leq -\inf_{g_1 \in O}l(g_0,g_1)\end{equation} 
We consider $\chi$ a smooth positive function equals to 0 outside O and equals to 1 on a small open ball centered in $g_1$ smaller than 1.  

We would like to apply the mechanism of Malliavin integration by parts to the measure
\begin{equation}h \rightarrow P_t[h\chi](g_0)\end{equation}
such that
\begin{equation}\vert P_t[\chi f^{(\alpha)}h](g_0)\vert \leq C t^{(-r_{(\alpha)})}\exp[{-l(g_0,g_1)+\delta\over t^{1/2k-1}}]\Vert h \Vert_\infty\end{equation}
for a small $\delta$. Since  (104) is true,  we have:
 
{\bf{
Theorem}}
When $t \rightarrow 0$
\begin{equation}\overline{Lim}_{t \rightarrow 0}t^{1/2k-1}Log \vert p_t(g_0,g_1)\vert  \leq- l(g_0,g_1) \end{equation}


\begin{thebibliography}{99}
\bibitem{1}S. Albeverio, R. Hoegh-Krohn: "Dirichlet forms and diffusion processes on rigged Hilbert spaces" {\it{Z.W. }}40, 1-57, 1977.
\bibitem{2}A. Asada: "Regularized calculus: an application of zeta regularization to infinite dimensional geometry". {\it{Int. J. Geom. Methods. Mod. Phys}},1. 107-157, 2004.
\bibitem{3}P. Auscher, P. Tchamitchian: {\it{Square root problem for divergence operators and related topics}}  {\it{Asterisque}} 249, Paris, France: S.M.F, 1998.
\bibitem{4}R. Azencott: "Grandes d\'eviations et applications". in {\it{Ecole de probabilit\'es de Saint-Flour}} {\it{P.L. Hennequin edt}}. {\it{Lect. Notes. Maths}} 778, Heidelberg, Germany: Springer, 2-176,1978.
\bibitem{5}F. Baudoin:{\it{An introduction to the geometry of stochastic flows}}. London, U.K.: Imperial College Press, 2000. 
 \bibitem{6}Y. Berezanskii: "The self-adjointes of elliptic operators with infinite number of variables". {\it{Ukrainian Math. Journal}} 27, 729-742, 1975.
\bibitem{7}J.M. Bismut: "Martingales, the Malliavin Calculus and hypoellipticity under general Hoermander's xondition".{\it{ Z.W.}} 56, 469-505, 1981.
\bibitem{8}J.M. Bismut: "Calcul des variations stochastiques et processus de sauts", {\it{Z. Wahr. Verw. Gebiete}}, 63, pp. 147-235, 1983.
\bibitem{9}J.M. Bismut: {\it{Large deviations and the Malliavin Calculus}} 45, Boston, U.S.A.: Birkhauser, 1984.
\bibitem{10}N. Bouleau, L. Denis: {\it{Dirichlet forms methods for Poisson measures and L\'evy processes with emphasis on the creation-annihilation techniques}}. {\it{Prob. Theory and stochastic modelling}} 76. Heidelberg, Germany: Springer, 2015.
\bibitem{11}K. Burdzy: "Some path properties of iterated Brownian motion" in {\it{Seminar on stochastic processes 1992}}(E. Cinlar and al eds). {\it{Progress Probab}} 33. Boston, U.S.A: Birkhauser. 67-87, 1993.
 \bibitem{12}J. Chazarain, A. Piriou: {\it{Introduction a la th\'eorie des \'equations au  d\'eriv\'ees partielles lin\'eaires}}. Paris, France:  Gauthier-Villars, 1981.
\bibitem{13}E.B. Davies:{\it{Heat kernels and spectral theory}}. Cambridge, England: Cambridge University Press, 1989.
\bibitem{14}E.B. Davies:"Uniformly elliptic operators with measurable coefficients", {\it{J. Funct. Anal.}}, 132, 141-169, 1995.
\bibitem{15}P.  Greiner: "An asymptotic expansion for the heat equation" {\it{Arch. Rat.Mech. Anal}} 41, 163-218, 1971.
\bibitem{16}T. Hida: {\it{Analysis of Brownian functionals}} {\it{Carleton Math. Lectures. Noters}} 13. Ottawa, Canada:  Carleton University Press, 1975.
\bibitem{17}L. Hoermander: "Fourier integral operators. I", {\it{Acta. Mathe.}}, 127, 78-183, 1971.
 \bibitem{18}L. Hoermander: {\it{The analysis of linear partial operators III}}.  Berlin, Germany:   Springer, 1984.
\bibitem{19}L. Hoermander: {\it{The analysis of linear partial operators IV}}. Berlin, Germany:  Springer, 1984.
\bibitem{20}N. Ikeda, S. Watanabe:{\it{Stochastic differential equations and diffusion processes}}., Amsterdam, Netheland: North-Hollan, 1989.
\bibitem{21}Y. Ishikawa : {\it{Stochastic Calculus of variations for jump processes}}, Basel, Schweiz: de Gruyter, 2012.
\bibitem{22}N. Jacob: {\it{Pseudo-differential operators and Markov Processes. I. Fourier Analysis and semigroups}}. London, England: Imperial College Press, 2001.
\bibitem{23}N. Jacob: {\it{Pseudo-differential operators and Markov processes. II. Generators and their potential theory}}. London, England: Imperial College Press, 2003.
\bibitem{24}N. Jacob: {\it{Pseudo-differential operators and Markov Processes. III. Markov processes and their applications}}. London, England: Imperial College Press, 2005.
\bibitem{25}N. Kumano-Go: "Phase space Feynman path integrals of higher order parabolic type with general functional as integrand". {\it{Bull. Sci. Math.}} 139, 495-537. 2015.
\bibitem{26}S. Kusuoka: "More recent theory of Malliavin Calculus" {\it{Sugaku Expositions}}5, 155-173, 1992.
\bibitem{27}P. Lax: "Asymptotic solutions of oscillatory initial value problems" {\it{Duke. Math. Journal}} 25, 627-646, 1957.
 \bibitem{28}R. L\'eandre:  "Extension du th\'eoreme de Hoemander a divers processus de sauts", PHD Thesis, Universit\'e de Besan\c con, France (1984)
\bibitem{29}R. L\'eandre: "Applications quantitatives et qualitatives du Calcul de Malliavin" In {\it{S\'eminaire Franco-Japonais}}(M. M\'etivier et S. Watanabe eds).
 {\it{Lectures Notes in Maths}}1322,Berlin, Germany: Springer, 109-134, 1988. English translation: {\it{Geometry of random motion}}(R. Durrett and M. Pinsky eds) {\it{Contemporary Maths}} 73. Providence, U.S.A. AMS. 173-196, 1988.
\bibitem{30}R. L\'eandre:  "Path integrals in noncommutative geometry" In "{\it{Encyclopedia of mathematical physics"}} 5J.P. Francoise and al eds). London, England: Elsevier, 8-12, 2006.
\bibitem{31}R. L\'eandre: "Malliavin Calculus of Bismut type without probability" {\it{Festchrift in honour of K. Sinha}} {\it{Proc. Indian. Acad.Sci: Math.Sci}}11-, 507-518, 2006.
\bibitem{32}R. L\'eandre: "Applications of the Malliavin Calculus without probability"{\it{WSEAS Trans Math}} 5, 1205-1206, 2006.
\bibitem{33}R. L\'eandre: "Malliavin Calculus of Bismut type in semi-group theory". {\it{Far East Jour. Math. Sciences}} 30, 1-26, 2008.
\bibitem{34}R. L\'eandre: "Wentzel-Freidlin estimates for an operator of order four". In {\it{2014 Int. Conp. Science. Comput. Int}} (B. Akhgar and al es). Los Alamitos, U.S.A.: IEEE Computer Society, 360-364, 2014 (IEEE X-plore).
\bibitem{35}R. L\'eandre: "Large deviation estimates for a non-markovian generator of big order" (E. Vagenas and al eds){\it{4th Int. Con. Math. Models. Physical. Sciences}}{\it{J. Phys: Conference series}} 633, paper 012085, 2015.
\bibitem{36}R. L\'eandre:  "Stochastic analysis for a non-markovian generator: an introduction", {\it{ Russian Journal of Mathematical Physics}}, 22, pp 39--52, 2015
\bibitem{37} R. L\'eandre: "The Ito-Stratonovitch for an operator of order four" (G. Budzhan an al eds) {\it{Festchrift in honour of A. Mukharjea, Ph. Feinsilver ans S. Mohammed}} {\it{Cotemporary Mathematics}} 668. Providence, U.S.A.: A.M.S.  165-169, 2016.
\bibitem{38}R. L\'eandre: "Large deviation estimates for an operator of order four with a potential"  in {\it{XII. Int. Work. Diff. Geometry and Application}} (R. Iordanescu edt){\it{Revue Roumaine de Math. Pures et Appliqu\'ees}} 61, 233-240, 2016.
\bibitem{39}R. L\'eandre: "Varadhan estimates for an operator of order four on a Lie group". In {\it{Control, decision and information technologies. 2016 }}(A. Vledma and al eds). IEEE-Xplore. Los-Alamitos, USA: IEEE.
\bibitem{40}R. L\'eandre: "A class of non-markovian generator of L\'evy type"in {\it{Pseudo-differential operators: groups, geometry and application}}(M.W. Wong and al eds). Berlin,Germany: Springer. 149-159. 2017. 
 \bibitem{41}R. L\'eandre, "Perturbation of the Malliavin Calculus of Bismut type of large order". in {\it{Phys. Mathematical ,Aspects of Symmetries}} (J.P. Gazeau and al eds). Heidelberg, Germany: Springer. 221-225, 2017.
\bibitem{42}R. L\'eandre, "Malliavin Calculus of Bismut type for an operator of order four on a Lie group". To appear {\it{Journal of Pseudo-differential operators and applications}}
 \bibitem{43}R. L\'eandre, "Bismut's way of the Malliavin Calculus of large order generators on a Lie group" In {\it{6th Int. Eurasian Conf. Math. Sciences and appl}}.( M. Tosun and al eds.) A.I.P.  1926, Melville, U.S.A. Ame. Inst. Phys.  020026.
\bibitem{44}P. Malliavin: "Stochastic calculus of variations and hypoelliptic operators"{\it{Pro. Symp. Stochastic Diff. equations}}. Tokyo, Japan. Kinokuyina. 195-263, 1978.
\bibitem{45}V.P. Maslov, M.V. Fedoriuk, {\it{Semiclassical approximation in quantum mechanics}}. Dordrecht, Netheland: Reidel, 1981.
\bibitem{46}S. Mazzucchi: "Infinite dimensional oscillatory integrals with polynomial phase and application to higher-order heat equation".To appear {\it{Potential analysis}}
\bibitem{47}S. Molchanov: "Diffusion processes and Riemannian geometry".Russ. Math.Surveys. 30, 1-63, 1975.
\bibitem{48} J. Neveu: {\it{Processus al\'eatoires Gaussiens}}{\it{S\'eminaire Maths Sup\'erieures}} 34.  Montr\'eal, Canada: Press. Univ. Montr\'eal, 1968.
\bibitem{49}D. Nualart:{\it{The Malliavin Calculus and related topics}} Berlin, Germany. Springer. 1995.
\bibitem{50}R.T. Seeley: "Complex powers of an elliptic operator". In {\it{Singular integrals}} {\it{Proc.Symp.Pure Mathematics}}. Providence, U.S.A.:A.M.S. 288-307.
\bibitem{51}A. Ter Elst, D.W. Robinson: "Subcoercive and subelliptic operators on Lie groups: variable coefficients" {\it{Pub. Inst. Math. Sciences}} 29, 745-801, 1993.
\bibitem{52}S.R.S. Varadhan: "Large deviations and applications". Philadelphia, U.SA.: S.I.A.M.  1984.
\bibitem{53}S. Watanabe: "Stochastic analysis and its applications" {\it{Sugaku Expositions}}5, 51-69, 1992.
\bibitem{54}A.D. Wentzel, M.J. Freidlin, {\it{Random perturbations of dynamical systems}}. Berlin, Germany: Springer, 1984.
 

 





 



\end{thebibliography}
\end{document}